\documentstyle[a4,12pt]{article}

\begin{document} 

\title{Q.H.I. spaces}
\author{V.Ferenczi}
\maketitle

\newcommand{\ti}[1]{\mbox{$\tilde{#1}$}}
\newcommand{\wti}[1]{\mbox{$\widetilde{#1}$}}
\newcommand{\T}{\mbox{$\tilde{T}$}}
\newcommand{\U}{\mbox{$\tilde{U}$}}
\newcommand{\nul}{\mbox{$\tilde{0}$}}
\newcommand{\al}{\mbox{$\alpha$}}
\newcommand{\be}{\mbox{$\beta$}}
\newcommand{\dd}{\mbox{$\delta$}}
\newcommand{\Zo}{\mbox{$Z^{\perp}$}}
\newcommand{\Yo}{\mbox{$Y^{\perp}$}}
\newcommand{\la}{\mbox{$\lambda$}}
\newcommand{\e}{\mbox{$\epsilon$}}
\newcommand{\norm}[1]{\mbox{$\|#1\|$}}
\newcommand{\YY}{\mbox{$Y'$}}
\newcommand{\TT}{\mbox{$T'$}}
\newcommand{\tutu}{\mbox{$t'$}}
\newcommand{\cvf}{\mbox{$\stackrel{w}{\rightarrow}$}}

\def\N{{\bf N}}

\abstract{A Banach space $X$ is said to be Q.H.I. if
eve\-ry 
infinite
 dimensional quo\-tient spa\-ce of $X$ is H.I.: that is,
 a space is Q.H.I. if the H.I. property is not only
stable passing to subspaces, but also passing
 to quotients and to the dual.
We show that Gowers-Maurey's space is Q.H.I.; then we
provide an example 
of a reflexive H.I. space ${\cal X}$
 whose dual is not H.I., from which it follows
that $\cal X$ is not Q.H.I..}

\section{Introduction}
In the following,
by {\em space (resp. subspace)}, we shall always mean
infinite dimensional Banach space (resp. closed subspace).
We recall the central definition of \cite{GM}
A Banach space is said to be
 {\em hereditarily indecomposable (or H.I.)} if
it has no decomposable subspace.
In others words, a space $X$ is H.I. if for any $\e>0$,
any subspaces $Y$ and $Z$ of $X$, there exist
vectors $y \in Y$, $z \in Z$, such that
$\norm{y}=\norm{z}=1$ and $\norm{y-z} \leq \e$.
 
We now give some notation that is useful for
the construction of Gowers-Maurey's space and similar
spaces.
Let $c_{00}$ be the space of sequences of scalars all but finitely many of
which are zero. Let $e_{1},e_{2},\ldots$ be its unit vector basis.
If $E\subset \N$, then we shall also use the letter
$E$ for the projection from $c_{00}$ to $c_{00}$ defined by
$E(\sum_{i=1}^\infty a_i e_i)=\sum_{i\in E}a_i e_i$. If
$E,F \subset \N$, then we write $E<F$ 
to mean that $\sup E<\inf F$.
 An {\em interval} of integers is a subset of
$\N$ of the form $\{a,a+1,\dots,b\}$ for some $a,b\in \N$.
For $N$ in $\N$, $E_{N}$ denotes the interval
$\{1,\ldots,N\}$. The {\em range} of a vector
 $x$ in $c_{00}$, written $ran(x)$, is the smallest interval $E$
such that $Ex=x$.
  We shall write $x<y$ to mean $ran(x)<ran(y)$; notice that
this is only defined on $c_{00}$. If
$x_1<\cdots<x_n$ we shall say that $x_1,\ldots,x_n$ are
 {\em successive}.

 Let $\cal X$ be the class of
 Banach sequence spaces
  such
that $(e_{i})_{i=1}^{\infty}$ is a norma\-lized
 bimonotone basis.
Notice that for $p \geq 1$,
$l_p$ is in $\cal X$.
 We denote by  $B(l_p)$ the unit ball
of $l_p \cap c_{00}$.

By a {\em block basis} in a space $X \in {\cal X}$
 we mean a sequence
$x_1,x_2,\dots$ of successive non-zero vectors in $X$
  (such a
sequence must be a basic sequence) and by a {\em block subspace}
 of a space
$X\in\cal X$ we mean a subspace generated by
 a block basis.

Let $f$ be the function $\log_{2}(x+1)$.
If 
 $X \in {\cal X}$, and all successive vectors $x_{1}, \ldots,
x_{n}$ in $X$ satisfy 
the inequality $f(n)^{-1} \sum_{i=1}^{n} \|x_{i}\| \leq
 \| \sum_{i=1}^{n} x_{i} \|$,
then we say that $X$ satisfies an {\em $f$-lower estimate}.
We denote by ${\cal X}(f)$
 the set of Banach spaces in $\cal X$
satisfying an
$f$-lower estimate.
 
A function $h:[1,+\infty) \rightarrow [1,+\infty)$
belongs to the Schlumprecht class $\cal F$ of functions
if it satisfies
the following five conditions:

(i) $h(1)=1$ and $h(x)<x$ for every $x>1$;

(ii) $h$ is strictly increasing and tends to infinity;

(iii) $\lim_{x\rightarrow\infty}x^{-q}h(x)=0$
 for every $q>0$;

(iv) the function $x/h(x)$ is concave and non-decreasing;

(v) $h(xy)\leq h(x)h(y)$ for every $x,y\geq 1$.

We notice that $f$ and $\sqrt{f}$ belong to ${\cal F}$.

Given $X$ in $\cal X$,
 given $g$ in $\cal F$, a functional $x^{*}$ in $X^*$
is an $(M,g)-form$ if
$\|x^{*}\|^* \leq 1$ and $x^{*}= \sum_{j=1}^{M} x_{j}^{*}$ for
some sequence $x_{1}^{*} < \cdots < x_{M}^{*}$ of successive
functionals such that $\|x_{j}^{*}\|^* \leq g(M)^{-1}$
 for each $j$.
Notice that if $X \in {\cal X}(f)$,
then every $(M,f)$-form is in the unit ball of $X^*$.

\

Let $X$ be a Banach space.
Let $y$ be a subspace of $X$.
We shall denote by $Id_Y$ the identity map from $Y$ to
$X$.
 An operator from $Y$ to $X$ is said to be {\em finitely singular} if
its restriction to some finite
 codimensional subspace is an isomorphism into.
It is said to be {\em infinitely singular} if it is not
finitely singular. An operator $S$ from $Y$ to $X$
 is said to be
 {\em strictly singular} if
the restriction of $S$ to a subspace is never
 an isomorphism
into. This is equivalent to saying that for any
$\e>0$, any $Z$, there exists $z$ in $Z$ such that
$\norm{S(z)} \leq \e \norm{z}$.

 By Proposition 2.c.4 of \cite{LT}, if $S$ is strictly
singular, then for every $Z$ and every $\e >0$, there
 exists
$Z' \subset Z$ such that $\norm{S_{/Z'}} \leq \e$.
 If $S$ is strictly singular and
$T$ is an operator, then $TS$ and $ST$ are strictly
singular whenever they are defined.
We denote by ${\cal S}(Y,X)$ the space of strictly
 singular operators
from $Y$ to $X$.

 Let $T$ be an operator from $Y$ to $X$.
By
Proposition 2.c.10 of \cite{LT}, if $T$ is of the form
$Id+S$, where $S$ is strictly singular,
 then it is finitely singular, that is,
an isomorphism on some finite codimensional subspace.
 It is said to be
an {\em $Id+S$-isomorphism} if it is  of
the form $Id+S$ and an isomorphism on the
whole subspace $Y$.
If $T$ is an $Id+S$-isomorphism, then so is
$T^{-1}=Id-(T-Id)T^{-1}$.
If $T$ and $U$ are $Id+S$-isomorphisms, then so is $TU$
when it is defined.
The subspaces $Y$ and $Z$ of $X$ are said to be
{\em $Id+S$-isomorphic} if there exists
 an $Id+S$ isomorphism
 from
$Y$ onto $Z$. 

We finally recall a classical notion:
two Banach spaces $X$ and $X'$ are 
{\em totally incomparable} if no subspace of $X$ is
 isomorphic
to a subspace of $X'$.

\section{Definition: Q.H.I. spaces}

\subsection{Definition}

A Ba\-nach spa\-ce $X$ is 
{\em quotient he\-re\-di\-ta\-ri\-ly inde\-com\-posa\-ble}
(or {\em Q.H.I.}) if eve\-ry 
infinite
 dimensional quo\-tient spa\-ce of $X$ is H.I..

\paragraph{Remark} Any Q.H.I. space is H.I. as a quotient space of
itself.

\paragraph{Proposition 1}
{\em Let $X$ be a Banach space. Assume that
for every infinite codimensional and infinite dimensional subspace $Y$,
$X/Y$ is  H.I.. Then $X$ is Q.H.I..}

\underline{Proof} It is enough to prove that $X$ is H.I.. Assume
$X$ is not H.I.. Then
$X$ contains a direct sum $Y \oplus Z$. Let $W$ be an infinite
dimensional and infinite codimensional subspace of $Y$
(for example the space generated by the even vectors of
 a basic sequence in $Y$). Then $X/W$ contains a space isomorphic to
the sum $Y/W \oplus Z$, so $X/W$ is not H.I..

\paragraph{Proposition 2}
{\em If $X$ is reflexive and Q.H.I., then so is $X^*$.}

\underline{Proof} A space $X$ is Q.H.I. if and only if any
subspace of a quotient of $X$ is indecomposable. For a reflexive space,
this property passes to the dual.

\subsection{Notation}
We shall often refer to the space defined by Gowers and
Maurey in \cite{GM}, and to the
 techniques developped in their article.
When we shall refer to a lemma of \cite{GM},
we shall write it  
 with the 
letters GM.

By {\em space} (resp. {\em subspace}), we shall always
mean infinite dimensional space (resp. subspace).
Let $X$ be a Banach space, $Y$ be a subspace of $X$.

An operator on $X$ is said to be {\em finitely singular} if
its restriction to some finite
 codimensional subspace is an isomorphism into.
It is said to be {\em infinitely singular} if it is not
finitely singular. It is said to be
{\em strictly singular} if all its restrictions are
infinitely singular.

An operator from $Y$ to $X$ is said to be an
 {\em $Id+S$-isomorphism} if it is an isomorphism of the
form $Id+S$, where $S$ is strictly singular.

The subspace $Y$  is said to be {\em quasi-maximal} if $Y$ and
any subspace $W$ of $X$ have $Id+S$-isomorphic subspaces.

We know from \cite{F} that if $X$ is H.I., then
 every subspace of $X$ is quasi-maximal; as an easy
consequence, if $X$ has
a quasi-maximal H.I. subspace then $X$ is H.I..

Two Banach spaces $X$ and $X'$ are 
{\em totally incomparable} if no subspace of $X$ is isomorphic
to a subspace of $X'$.

For $p \geq 1$, $B(l_p)$ denotes the unit ball
of $l_p \cap c_{00}$. An interval of integers $E$
is a subset of $\N$ of the form $\{m,m+1,\ldots,n\}$.
The letter $E$ will also denote the projection defined on
$c_{00}$ by $E(\sum_{i=1}^{+\infty} \la_i e_i)=\sum_{i \in E}\la_i e_i$.
The {\em range} of a vector $x$ in $c_{00}$, denoted $ran(x)$
 is the smallest interval
$E$ such that $Ex=x$. We say that $E$ and $F$ are {\em successive} if
$\max E<\min F$. Two vectors in $c_{00}$ are said
to be {\em successive} if their ranges are successive.

 Let
${\cal X}$ be the set of Banach spaces
with a bimonotone basis. 
Let $f$ be the function $x \rightarrow log_2(x+1)$.
A space $X$ in ${\cal X}$
 is said to {\em satisfy an $f$-lower estimate} if
any successive vectors $x_1<\cdots<x_n$ in $X$ satisfy the inequality
$\norm{\sum_{i=1}^n x_i} \geq f(n)^{-1}\sum_{i=1}^n \norm{x_i}$.
Let ${\cal X}_f$ be the set of Banach spaces
satisfying an
$f$-lower estimate.
An $(M,f)$-form is a vector of $X^*$ of the form
$f(M)^{-1}\sum_{i=1}^M x_i^*$, where the $x_i^*$'s are
successive vectors in the unit ball of $X^*$. If $X \in {\cal X}_f$,
then every $(M,f)$-form is in the unit ball of $X^*$.

\section{Successive vectors in a quotient space}

In this section, we want to define a notion of successive vectors and
of $l_1^{n+}$-vectors in
a quotient space of a reflexive space in ${\cal X}_f$. As a quotient space needs not have a basis, we cannot
use the usual definition.
Furthermore, our new definition does not necessarily coincide
with the usual one when the quotient space has a basis.
The first lemmas are rather technical; they allow us to
show Lemma 7, which is, with the definition of successive vectors and
of $l_1^{n+}$-averages in a quotient space,
 the main tool used afterwards.

\subsection{Notation}

Let $0<\Delta<1$.
\paragraph{Definition}
Let $W$ be a Banach space. Let $(w_n)$ and $(w'_n)$ be two se\-quen\-ces
 in $W$.
They are said to be $\Delta$-{\em equi\-va\-lent} if for all $n$,
 \[\norm{w_n-w'_n} \leq \Delta^{n} \inf(\norm{w_n},\norm{w'_n}).\]

\paragraph{Definitions}
Let $W$ be a Banach space with a basis.
A sequence $(w_n)$ in $W$ is said to be
{\em almost successive} if it is $\Delta$-equivalent to a sequence
of successive vectors in $W$.
Let $V$ be a subset of $W$.
A sequence $(w_n)$ in $W$ is said to be
{\em almost in $V$} if it is $\Delta$-equivalent to a sequence
of vectors in $V$.

\paragraph{Lemma 1}

{\em  Let $W$ be a Banach space in ${\cal X}$.
Let $(w_i)$ be a successive sequence in $W$
and $(w'_i)$ be $\Delta$-equivalent to $(w_i)$.
 Then \[\norm{\sum_{i=1}^n w'_i} \leq
 (1+2\Delta)\norm{\sum_{i=1}^n w_i}.\]}

\underline{Proof}
By definition,

\[ \norm{\sum_{i=1}^n w'_i} \leq
 \norm{\sum_{i=1}^n w_i}+\sum_{i=1}^n \Delta^{i} \norm{w_i}, \]
so, as the basis in $W$ is bimonotone,
\[\norm{\sum_{i=1}^n w'_i} \leq (1+2\Delta)
\norm{\sum_{i=1}^n w_i}.\]

\paragraph{Lemma 2}
{\em Let $W$ be a space in ${\cal X}$. Let $(w_n)$ be
 a sequence in $W$ such that $w_n/\norm{w_n} \cvf 0$. Then $(w_n)$ has
an almost successive subsequence.}

\underline{Proof}
We may assume $(w_n)$ is a norm $1$ sequence. 
Assume we have already selected $w_{n_1},\ldots,w_{n_{k-1}}$
satisfying the conclusion; let
$t_1,\ldots,t_{k-1}$ be the associated successive sequence.
Let $E$ be an interval containing the first vector of the basis
 and the range of
$t_1+\cdots+t_{k-1}$. There exists $n_k$ such that
$\norm{Ew_{n_k}} \leq {\Delta}^{k}/4$.
 Let $t'_{k}=w_{n_k}-Ew_{n_k}$.
 There exists
an interval $F$ such that $Ft'_{k}$ is equal
 to $t'_k$ up to
 ${\Delta}^{k}/4$.
If we let $t_{k}=Ft'_{k}$, we have that
$t_{k}>t_{k-1}$,
and $\norm{t_{k}-w_{n_k}} \leq {\Delta}^{k}/2$, from
 which it follows that
\[\norm{t_{k}-w_{n_k}} \leq
 {\Delta}^k\inf(\norm{t_k},\norm{w_{n_k}}),\]
so the se\-quen\-ce $w_{n_1},\ldots,w_{n_k}$ 
and the suc\-ces\-sive se\-quence $t_1,\ldots,t_k$
are $\Delta$-equi\-va\-lent.
\

>From now on, $X$  stands for a reflexive
Banach space in ${\cal X}_f$,  $Y$ for an infinite dimensional and
infinite codimensional subspace of $X$, and
we assume that $\Delta=10^{-4}$.

\subsection{Successive vectors in $X/Y$}

\paragraph{Definitions}

Given $z$ in $X$, we denote by $\hat{z}$ its class in $X/Y$.

A {\em couple} denotes a sequence $((\wti{z_n},z_n^*))$ in 
$X \times X^*$.

For an element $(\wti{z},z^*)$ in
$X \times X^*$, the {\em range} is the smallest interval
containing $ran(\wti{z})$ and $ran(z^*)$.

A couple is said to be {\em successive} if the sequence of ranges
is successive.

Let $(z_n)_{n \in \N}$ be a sequence in $X/Y$.
A couple $(\wti{z_n},z_n^*)$ is said to be
{\em associated to $(z_n)$} if 
$(\widehat{\wti{z_n}})$ and $(z_n)$ are $\Delta$-equivalent
and $(z_n^*)$ is almost in $B(\Yo)$.

Let $\la >1$. A couple associated to $(z_n)$ is said to be
{\em $\la$-associated to $(z_n)$} if for all $n$,
$\norm{\wti{z_n}} \leq \la \norm{z_n}$ and
$z_n^*(\wti{z_n}) \geq (1/3) \norm{z_n}$.

\paragraph{Definition}
Let $\la>1$. Let $(z_n)_{n \in \N}$ be a sequence in $X/Y$. We say
that $(z_n)$ is {\em $\la$-successive} if it is 
$\la$-associated to some successive couple  $((\wti{z_n},z_n^*))$.
We shall then call $(\wti{z_n})$ a $\la$-lifting,
$(z_n^*)$ a $\la$-functional.

\paragraph{Remark} If $(z_n)$ is $\la$-successive, then for
any $\la' \geq \la$, any subsequence $(z'_n)$ of
$(z_n)$, $z'_n$ is $\la'$-successive.

\paragraph{Lemma 3}
{\em Let $(z_n)$ be a $400$-successive
 sequence in $X/Y$. Then for all $n$,
\[ (4f(n))^{-1}\sum_{i=1}^n \norm{z_i} \leq \norm{\sum_{i=1}^n z_i}.\]}

\underline{Proof} Let $(z_n)$ be a $400$-successive sequence.
Let $(z_n^*)$ in $X^*$ be the successive functional
 associated to $(z_n)$,
 let $(y_n^*)$ in $B(Y^{\perp})$ be the almost successive functional 
$\Delta$-equivalent to $(z_n^*)$.
By Lemma 1 , we have that 
 \[\norm{\sum_{i=1}^n y_i^*} \leq (1+2\Delta)\norm{\sum_{i=1}^n z_i^*},\]
so
 \[\norm{\sum_{i=1}^n y_i^*} \leq (1+2\Delta)(1+\Delta)f(n) \leq 
(1+4\Delta)f(n),\]

It follows that
\[ \norm{\sum_{i=1}^n z_i}
 \geq ((1+4\Delta)f(n))^{-1}\left(\sum_{i=1}^n y_i^*\right)
\left(\sum_{i=1}^n z_i\right),\] 
\[ \norm{\sum_{i=1}^n z_i} \geq ((1+4\Delta)f(n))^{-1}
\left((\sum_{i=1}^n z_i^*)(\sum_{i=1}^n \wti{z_i})-
(\sum_{i=1}^n\norm{z_i^*-y_i^*})(\sum_{i=1}^n \norm{\wti{z_i}})\right),\]

\[ \norm{\sum_{i=1}^n z_i} \geq ((1+4\Delta)f(n))^{-1}
(1/3-800\Delta)\sum_{i=1}^n \norm{z_i} 
\geq (4f(n))^{-1}\sum_{i=1}^n \norm{z_i}.\]

\paragraph{Lemma 4}
{\em Let $(z_n)$ be a sequence in $X/Y$ such that
$z_n/\norm{z_n} \cvf 0$. Then $(z_n)$ has a $3$-successive subsequence.}

\underline{Proof}
We may assume that $(z_n)$ is a norm $1$ sequence.
Let $(\wti{z_n}')$ be a lifting for $(z_n)$
such that $\norm{\wti{z_n}'} \leq 1+{\Delta}^{n}$.
The sequence $(\wti{z_n}')$ is bounded, so, passing to
a subsequence, we may assume that $(\wti{z_n}')$ converges
weakly. Let $y$ be the weak limit of $(\wti{z_n}')$.
The vector $y$ has norm $1$, and belongs to $Y$, because for every
$y^*$ in $\Yo$, $y^*(\wti{z_n}')=y^*(z_n)$ tends to $0$.

Let $\wti{y_n}=\wti{z_n}'-y$. We have 
\[ \norm{\wti{y_n}} \leq 2+{\Delta}^{n}.\]
We have  that $\wti{y_n} \cvf 0$ so
passing to a further subsequence, we may assume by Lemma 1
that $(\wti{y_n})$ is almost successive; let 
$(\wti{z_n})$ be a successive sequence equi\-valent to
$(\wti{y_n})$. We have that
\[ \norm{\wti{z_n}} \leq (1+{\Delta}^{n}) \norm{\wti{y_n}}
\leq 3\norm{z_n}.\]

Now let $(y_n^{\prime *})$ be a dual sequence in $B(Y^{\perp})$
such that for all $n$,
$y_n^{\prime *}(z_n)=1$. Passing to a subsequence, we may assume that
$y_n^{\prime *} \cvf y^*$ and that for all $n$,
$|y^*(z_n)| \leq 1/6$. Let $y_n^*=1/2(y_n^{\prime *}-y^*)$.
 We have that $y_n^* \in B(\Yo)$ and
\[y_n^*(z_n)=1/2(y_n^{\prime *}(z_n)-y^*(z_n)) \geq 5/12,\]
for $n$ greater than some $N$. Furthermore, $y_n^* \cvf 0$.
Passing to a new subsequence, by Lemma 3, we may assume that
$(y_n^*)$ is almost successive. Let $(z_n^*)$ be the
equivalent successive sequence. Then
\[ z_n^*(\wti{z_n}) \geq y_n^*(z_n)-\norm{z_n^*-y_n^*}\norm{\wti{z_n}},\]
so that
\[ z_n^*(\wti{z_n}) \geq 5/12-3.\Delta^{n} \geq 1/3.\]

\paragraph{Definition}
Let $\la>1$. Let $(z_n)$ be a sequence in $X/Y$. For any finite
subset $E$ of $\N$,
$z_E$ denotes the sum $\sum_{i \in E}z_i$.
For $j \geq 0$, $n \geq 0$, 
let $E_n(j)$ be the interval $[2^jn+1,2^jn+2^j]$.
We say that $(z_n)$ is {\em $\la$-supersuccessive}
if for every $j$, the sequence $(z_{E_n(j)})_{n \geq 1}$
 is $\la$-successive.

\paragraph{Lemma 5}
{\em Let $(z_n)$ be a norm $1$ sequence in $X/Y$ such that
$z_n \cvf 0$. Then $(z_n)$ has a
$3$-supersuccessive subsequence.}

\underline{Proof}
It  is enough to prove by induction that there exists an
inclusion decreasing
sequence of subsequences
 $(z_n^1)_{n \in N},\ldots,(z_n^k)_{n \in \N}$ such that
for all $k$, all $j\leq k$, the sequence
$(z^k_{E_n(j)})_{n \geq 1}$ is $3$-successive.
 Indeed, the result then follows by taking a diagonal
subsequence.

Assume the induction hypothesis is true for $k-1$.
The sequence $(z^{k-1}_{E_n(k)})$ converges weakly
 to $0$ and by Lemma 3,
 it is bounded 
below, so by Lemma 4,
it admits a $3$~-~successive subsequence 
$(z^{k-1}_{E_{n_i}(k)})_{i \in \N}$.
Let $E=\cup_{i} E_{n_i}$.
We let $(z_n^k)_{n \in \N}$ be the sequence
$(z_n^{k-1})_{n \in E}$.
 It is a subsequence of $(z_n^{k-1})_{n \in \N}$.
Furthermore, it is easy to check that
the sequence $(z^k_{E_n(j)})_{n \geq 1}$ is
still successive for $j \leq k-1$, and by construction,
this is also true for $j=k$.

\subsection{$l_1^{n+}$ vectors in $X/Y$}
\paragraph{Definitions: $l_1^{n+}$-vectors}
\
Let $\e_0=1/40$, let $C>0$.

An {\em $l_1^{n+}$-vector with constant $C$ in $X$} is
a vector $x$ of the form $\sum_{i=1}^n x_i$ such that
the sequence $(x_i)$ is successive and
for all $i$, $\norm{x_i} \leq C \norm{x}/n$.

An {\em $l_1^{n+}$-average in $X$} is a norm $1$
$l_1^{n+}$-vector in $X$.

An {\em $l_1^{n+}$-vector in $X/Y$} is
a vector $x$ of the form $\sum_{i=1}^n x_i$ such that
the sequence $(x_i)$ is $3$-successive and
for all $i$, $\norm{x_i} \leq (1+\e_0) \norm{x}/n$.

Let $x$ be a
 $l_1^{n+}$-vector in $X/Y$ of the form $\sum_{i=1}^n x_i$.
An {\em $l_1^{n+}$-lifting} for $x$ is a lifting
$\ti{x}=\sum_{i=1}^n \wti{x_i}$, whe\-re 
$(\wti{x_i})_{1 \leq i \leq n}$ is a $3$-lifting
 for $(x_i)_{1 \leq i \leq n}$.

An {\em $l_1^{n+}$-ave\-ra\-ge in $X/Y$}
is an $l_1^{n+}$-vec\-tor in $X/Y$
the $l_1^{n+}$-lifting of which
 has norm $1$.

\paragraph{Lemma 6}
{\em Let $n \in \N$.
 Let $x=\sum_{i=1}^n x_i$ be a $l_1^{n+}$-average in $X/Y$,
let $\ti{x}$ be an $l_1^{n+}$-lifting for $x$. Then
 $\ti{x}$ is a $l_1^{n+}$-average with constant
$4$ and $\norm{\ti{x}} \leq 4\norm{x}$.}

\underline{Proof}
We know that the vectors 
$(\wti{x_i})$ are successive in $X$.
Now

\[ \norm{x-\widehat{\ti{x}}} \leq
 \sum_{i=1}^{n} \norm{x_i-\widehat{\wti{x_i}}}
 \leq \sum_{i=1}^{n} {\Delta}^{i} \norm{x_i}
\leq 2\Delta(1+\e_0)/n \norm{x}.\]
It follows that 
\[\norm{x} \leq (1-2\Delta (1+\e_0)/n)^{-1}\norm{\widehat{\ti{x}}}
\leq (1-2\Delta (1+\e_0)/n)^{-1}\norm{\ti{x}}.\]
So \[\norm{\wti{x_i}} 
\leq 3\norm{x_i}\leq 3(1+\e_0)\norm{x}/n
\leq 4\norm{\ti{x}}/n,\] 
so $\ti{x}$ is a $l_1^{n+}$ with constant
$4$. Furthermore,

\[\norm{\ti{x}} \leq \sum_{i=1}^n \norm{\wti{x_i}} \leq
3\sum_{i=1}^n \norm{x_i} \leq 3(1+\e_0)\norm{x} \leq 4\norm{x}.\]

\paragraph{Lemma 7}
{\em Let $Z \subset X/Y$. Then $Z$ contains
a $4$-successive basic sequence $(Z_n)$ of $l_1^{2^n +}$-averages, such
that for all $n$, the $4$-lifting $\wti{Z_n}$ is 
a $l_1^{2^n +}$-lifting. Such a sequence will be called
an {\em average basic sequence in $Z$}.}

\underline{Proof}
As $Z$ is reflexive, there exists a basic sequence $(z_n)$ of unit
vectors such that $z_n \cvf 0$. By Lemma 5,
we may assume that $(z_n)$ is $3$-supersuccessive.

We now prove the result by induction. Assume we have 
chosen the first $n-1$ terms  $Z_1,\ldots,Z_{n-1}$
 of the sequence.
Let $(\wti{Z_j},Z_j^*)_{1 \leq j \leq n-1}$ be the associated
$4$-successive couple.
Let $T$ satisfy 
$4 f(2^{nT}) < (1+\e_0)^{T}$ and let $N=nT$.
 Let $Z'_1,\ldots,Z'_{2^N}$ be the sequence
$Z_{2^Nr+1},\ldots,Z_{2^Nr+2^N}$, for some $r$ such that,
with this notation, 
 \[ \min_{j \leq n, p \leq 2^{N-j}}
      ran((\wti{Z'_{E_p(j)}},{Z'_{E_p(j)}}^*)) > 
      ran((\wti{Z_{n-1}},Z_{n-1}^*)).\]

For any $j \leq N$,
the sequence $(Z'_{E_n{j}})$ is $3$-successive.
Now assume no decomposition of any $Z'_{E_p(j)}$ as 
$\sum_{k=1}^{2^n} Z'_{E_k(j-1)}$ is an $l_1^{2^n +}$
 decomposition. It follows by induction that
\[ \norm{Z'_{E_p(j)}} \leq (2^n/(1+\e_0))^j,\]
so that
\[ \norm{Z'_{E_1(T)}} \leq 2^N/(1+\e_0)^{T}.\]
But $(z_n)$ is $3$-successive so by Lemma 3,
\[ \norm{Z'_{E_1(T)}}=
\norm{Z'_{1}+\cdots+Z'_{2^N}} \geq 2^N/4f(2^N),\]
a contradiction by choice of $N$.

It follows that one of the $Z'_{E_p(j)}$ is a 
$l_1^{2^n+}$-vector; we choose $Z_n$ to be the associated
$l_1^{2^n+}$-average. By choice of
$r$, the couple $(\wti{Z_n},Z_n^*)$
 satisfies $(\wti{Z_n},Z_n^*) > (\wti{Z_{n-1}},Z_{n-1}^*)$. 

\section{Rapidly Increasing Sequences}

\subsection{R.I.S.-vectors in a space of the class ${\cal X}_f$}
 We now define R.I.S.-vectors in a Banach space $X$ in ${\cal X}_f$. 
In fact, the properties of R.I.S. are not interesting in all
spaces in ${\cal X}_f$, but only on those of Gowers-Maurey's type;
we give a sense to this expression in the following paragraph,
defining GM-type spaces, and
then state several lemmas true in GM-type spaces.

\paragraph{Definitions}

Let $J$ be a set of integers $\{j_n,n \in \N\}$, such that
$f(j_1)>256$ and for all $n$, $\log \log \log j_{n+1} \geq 4 j_n^2$, let
$K=\{j_1,j_3,j_5,\ldots\}$, let $L=\{j_2,j_4,j_6,\ldots\}$.

Let $L' \subset L$. 
An {\em $L'$-sequence} is a successive sequence
$x_1^*<\cdots<x_k^*$ with $k \in K$, such that for all $i$, $x_i^*$ is
an $(M_i,f)$-form where $M_i$ is an element in $L'$ greater than $j_{2k}$.
An {\em $L'$-sum} is a vector of the form
$1/\sqrt{f(k)}\sum_{i=1}^k x_i^*$, where $x_1^*,\ldots,x_k^*$
 is a $L'$-sequence.

A space $X$ is {\em of GM-type} if it belongs to ${\cal X}_f$,
 and if there exists a set $\cal S$ of $L$-sums in $B(X^*)$ such
 that every vector in $X$
has either the supremum norm, or is normed by an $(M,f)$-form, or
  by an element in $\cal S$.

\

Notice that Gowers-Maurey's space is of GM-type (the set $\cal S$
being the set of special sums), but that a GM-type space needs not be H.I..
Notice also that a space of GM-type satisfies a $f$-lower estimate,
 so by Lemma GM3, every block subspace contains
$l_1^{n+}$-vectors with arbitrary constants and lengths.

\paragraph{Definition}
We recall that a {\em R.I.S. of length $N$ with constant $\e$ in $X$}
is a successive sequence $(x_i)_{i=1}^N$
of $l_1^{n_i}$-averages with constant $C$ in $X$ 
 such that
$(n_i)$ satisfies
$n_1 \geq
4(1+\e)M_f(N/\e')/\e'$ and
 $\e'/2 f(n_k)^{1/2} \geq |ran(x_{k-1})|$ for $k=2,\ldots,N$,
where $\e'=\min\{\e,1\}$ and $M_f(x)=f^{-1}(36x^2)$.

A R.I.S.-vector is a non-zero multiple of the sum of a R.I.S..

We now show some lemmas very similar to those of \cite{GM};
we have to state them because we shall use different constants,
and because they can be applied to any GM-type space, which 
will be useful in the last part of the article.

\paragraph{Lemma 8}
{\em Let $X$ have GM-type. Let $N \in L$, let
$n \in [\log N,\exp N]$, let $(x_i)_{i=1}^N$ 
be a R.I.S. of length $M$ with constant $\e$ in $X$.
Then $\norm{\sum_{i=1}^n x_i} \leq (1+\e+\e')nf(n)^{-1}$.}

\underline{Proof}
As $X$ has GM-type, all the hypothesis of Lemma GM7 are
satisfied. The conclusion then follows from Lemma GM9.

\paragraph{Lemma 9}
{\em Let $X$ have GM-type. Let $N \in L$. Let $M=N^{1/40}$.
Let $x_1,\ldots,x_N$ be a R.I.S. in $X$ with constant $4$. Then
$\sum_{i=1}^N x_i$ is a $l_{1+}^M$-vector with constant $6$.}

\underline{Proof}
Let $m=N/M$, let $x=\sum_{i=1}^N x_i$, and for $1 \leq j \leq M$
let $y_j=\sum_{i=(j-1)m+1}^{jm}x_i$. By Lemma 8,
$\norm{y_j} \leq 5mf(m)^{-1}$, while
$\norm{\sum_{j=1}^m y_j}=\norm{x} \geq Nf(N)^{-1}$, so
$x$ is an $l_1^{M+}$-vector with constant
$5f(N)/f(m) \leq 6$.

\paragraph{Lemma 10}
{\em Every GM-type space is reflexive.}

\underline{Proof}
We follow the proof 
that Gowers-Maurey's space is reflexive. Let
$X$ have GM-type. We show that the canonical basis $e_1,e_2,\ldots$
of $X$
is boundedly complete and shrinking (we refer to \cite{LT} 
for the definition of these notions). It follows from
the fact that $X$ belongs to ${\cal X}_f$ that
the basis is boundedly complete. Now assume it is not
shrinking. Then we can find $\e>0$, a norm-$1$ functional
$x^*$ in $X^*$, and a sequence of successive vectors
$x_1,x_2,\ldots$ such that $x^*(x_n) \geq \e$ for every
$n$. It follows that $\sum_{n \in A}x_n$ is a $l_{1+}^{|A|}$-vector
with constant $\e^{-1}$ for any $A \subset \N$.
Given $N \in L$, we may construct a R.I.S. $y_1,\ldots,y_N$ with constant
$\e^{-1}$ with such $l_1^{n+}$-vectors, and we have
$x^*(y_1+\cdots+y_N) \geq \e N$. For $N$ large enough, this
contradicts Lemma~8. 

\paragraph{Definition} Let $X$ have GM-type.
Let $x_1^*,\ldots,x_k^*$ 
be an $L$-sequence 
of lenght $k$, where each $x_i^*$
 is an $(M_i,f)$-form. 
A  sequence of successive vectors $x_1<\cdots<x_k$ in $X$
 is said to be
{\em a R.I.S. associated to $x_1^*,\ldots,x_k^*$}
 if for every $i$,
$x_i$ is a normalized R.I.S. of lenght $M_i$ and constant
$4$, $M_1=j_{2k}$ 
and $1/2 f((M_i)^{1/40})^{1/2} \geq |ran(x_{i-1})|$.

\paragraph{Remark}
Because of the increasing condition, and by Lemma 9,
 a R.I.S. associated to an $L$-sequence of length $k$ is
a R.I.S. with constant $6$.

\paragraph{Lemma 11}
{\em Let $X$ have GM-type.
Let $x$ be a norm $1$ R.I.S.-vector in $X$ of length $N_1 \in L$ and constant
$4$ and
let $x^*$ be an $(N_2,f)-form$ in $X^*$ with $N_2 \in L$,
 and assume
$N_1 \neq N_2$. Let $k \in K$ be such that
 $N_1 \geq j_{2k}$, $N_2 \geq j_{2k}$. 
Then for every interval $E$, $|x^*(Ex)| \leq 1/k^2$.}

\underline{Proof} 
The proof relies on
 Lemmas GM4 and GM5, which
we may apply since $X$ satisfies an $f$-lower estimate.
First, by Lemma 9, $x$ is a
$l_1^{N'_1 +}$-average with constant $6$, where $N'_1=N_1^{1/40}$.

Now if $N_2<N_1$, then $N_2<N'_1$ by the lacunarity of $L$. By Lemma
GM4, 
$ |x^*(Ex)|=|(Ex^*)(x)| \leq 18/f(N_2) \leq 18/f(j_{2k}) \leq k^{-2}$.

If $N_1<N_2$, then $M_f(N_1)<N_2$ by the lacunarity of $L$,
so   we may apply Lemma GM5 to $x'$, the sum of the R.I.S. 
whose
normalized sum is $x$. It follows that $|x^*(Ex')| \leq 5$, while
$\norm{x'} \geq N_1/f(N_1) \geq 5k^{2}$.

\paragraph{Lemma 12}
{\em Let $X$ have GM-type. Let $k \in K$.
Let $x_1<\ldots<x_k$ in $X$
 be a R.I.S. associated to some
$L$-sequence.
Let $x=\sum_{i=1}^k x_i$.

Assume that for every $L$-sum
in $\cal S$, every interval $E$,
 $|z^*(Ex)| \leq 1/4$. Then
\[ \norm{x} \leq 7k/f(k). \]}

\underline{Proof}
By  Lemma GM9 of \cite{GM},
for any $K_0 \subset K$,
there exists a function $g_{K_0}$ belonging to the "Schlumprecht" class
$\cal F$ of functions, such that
$g_{K_0}(x)=\sqrt{f(x)}$ when $x \in K_0$, and
$g_{K_0}(x)=f(x)$ when $x \in [\log N,\exp N]$, for any
$N \in J \setminus K_0$.

Now take $K_0=K \setminus \{k\}$, and let $g$ be the function
associated to $K_0$. By the hypothesis,
for any interval $E$,
\[ \norm{Ex} \leq 1 \vee \sup\{|x^*(Ex)|:M \geq 2,x^*\ (M,g)-form\}.\]
As $x_1,\ldots,x_k$ is a R.I.S. with constant $6$, it is an easy
consequence of Lemma GM7 that $\norm{x} \leq 7k/f(k)$.

\subsection{R.I.S.-vectors in a quotient space}
We finally define R.I.S. in a quotient space $X/Y$ where $X$ has GM-type
and link them with
the R.I.S. in $X$.

\paragraph{Definition} Let $X$ have GM-type. Let $Z \subset X/Y$.
Let $(Z_n)$ be an average basic sequence in $Z$.

A {\em R.I.S. of length $M$ in $Z$}
is a subsequence $(Z_{n_i})_{i=1}^N$ of $(Z_n)$
such that $2^{n_1} \geq 16M_f(N)$ and
$1/2 f(2^{n_i})^{1/2} \geq ran(\wti{Z_{n_{i-1}}})$ for $i=2,\ldots,N$.

\paragraph{Lemma 13}

{\em Let $X$ have GM-type. Let $Z \subset X/Y$. Let $M \in L$. Let
$z=\sum_{i=1}^M Z_{n_i}$ be a R.I.S.-vector of length $M$ in $Z$.
 Let
 $\ti{z}=\sum_{i=1}^M \wti{Z_{n_i}}$
 be the sum of the
 sequence of the associated
$l_1^{n_i +}$-liftings. Then 
$\ti{z}$ is a R.I.S.-vector in $X$ with constant $4$.
Furthermore $\norm{\ti{z}} \leq 80 \norm{z}$.}

\underline{Proof} The first part is a direct consequence
of Lemma 6 and of the definitions. For the second, notice that

\[\norm{\ti{z}}=\norm{\wti{Z_{n_1}}+\cdots+\wti{Z_{n_M}}}
\leq 5M/f(M)\]
 by Lemma 8, while
\[ \norm{z}=\norm{Z_{n_1}+\cdots+Z_{n_M}} 
\geq (\sum_{i=1}^M \norm{Z_{n_i}})/4f(M)\]
by Lemma 3.
As $\norm{Z_{n_i}} \geq(1/4) \norm{\wti{Z_{n_i}}} \geq 1/4$ by Lemma 6,
it follows that
$\norm{\ti{z}} \leq 80\norm{z}$.

\section{Gowers-Maurey's space is Q.H.I.}
We refer to \cite{GM} for the definition of Gowers-Maurey's space,
which we shall denote by $X$. Gowers and Maurey have proved that $X$ 
is H.I.. Now let $Y$ be any infinite dimensional and infinite 
codimensional subspace of $X$. We shall prove that $X/Y$ is H.I..
We first show a Lemma similar to Lemma GM2.

\paragraph{Lemma 14}
{\em Let $x_1^*,\ldots,x_k^*$ be a special sequence in $X$.
Let $x_1<\ldots<x_k$ be a R.I.S associated to $x_1^*,\ldots,x_k^*$.
Let $x=\sum_{i=1}^k x_i$.

Assume that for
every interval $E$,
 $|(\sum_{i=1}^k x_i^*)(Ex)| \leq 2$, then
\[ \norm{x} \leq 7k/f(k). \]}

\underline{Proof}
We already know that $X$ has GM-type. By Lemma 12, it is enough
to prove that for any special function $z^*$, every interval $E$,
$|z^*(Ex)| \leq 1/4$. We follow the proof of Lemma GM2.
Let $z^*$ be such a function of the form
$f(k)^{-1/2}\sum_{i=1}^k z_i^*$. Let $t$ be maximal
such that $z_t^*=x_t^*$, or $0$  if no such $t$ exists.
Assume $i \neq j$ or one of $i,j$ is greater than $t+1$.
Then since $\sigma$ is an injection, we can find
$l \neq l'$ in $L$ such that
$z_i^*$ is an $(l,f)$-form and $x_j$ is a norm $1$ R.I.S.-vector
of length $l'$. It follows then from Lemma~11 that
$|z_i^*(Ex_j)| \leq k^{-2}$.
 Now choose an interval $F$ such that
\[\left|\left(\sum_{i=1}^t z_i^*\right)(Ex)\right|=
\left|\left(\sum_{i=1}^k x_i^*\right)(Fx)\right|
\leq 2.\]
It follows that
\[\left|\left(\sum_{i=1}^k z_i^*\right)(Ex)\right|
\leq 2+|z^*_{t+1}(x_{t+1})|+k^2.k^{-2} \leq 4,\]
and that
$|z^*(Ex)| \leq 4f(k)^{-1/2} < 1/4$.

\paragraph{Conclusion}

Let $Z$ and $Z'$ be two subspaces of $X/Y$. We want to
prove that their sum is not direct.
Lemma 7 allows us to consider an average basic sequence $(Z_n)$
(resp. $(Z'_n)$)
  in $Z$ (resp. $Z'$).
Let $\delta>0$, let $k \in K$ be 
such that $2.10^4/\sqrt{f(k)}\leq \delta$.

Let $M_1=j_{2k}$. We may build with vectors of the average
basic sequence $(Z_n)$  a R.I.S. vector
in $Z$ of length $M_1$, of the form
$z_1=\sum_{j=1}^{M_1}z_1^j$.
Let $\wti{z_1}=\sum_{j=1}^{M_1} \wti{z_1^{j}}$.
By Lemma 13, $\wti{z_1}$ is the sum of a R.I.S. with constant
 $4$, and we may assume it is of norm $1$. 

Let $(z_1^{j*})$ be the successive functional in $X^*$
associated to $(z_1^j)$,
$(y_1^{j*})$ be the almost successive functional in $B(\Yo)$ 
$\Delta$-equivalent to
$(z_1^{j*})$. Let $z_1^{\prime*}=f(M_1)^{-1}(\sum_{j=1}^{M_1}z_1^{j*})$. We
 have that
$(1+\Delta)^{-1}z_1^{\prime*}$  is a $(M_1,f)$-form, and
\[z_1^{\prime*}(\wti{z_1})=f(M_1)^{-1}\sum_{j=1}^{M_1}z_1^{j*}(\wti{z_1^j}),\]
so \[z_1^{\prime*}(\wti{z_1})\geq (3f(M_1))^{-1}\sum_{j=1}^{M_1}\norm{z_1^j},\]
and by Lemma 13 and Lemma 2,
\[z_1^{\prime*}(\wti{z_1}) \geq 1/240 \sum_{j=1}^{M_1}\norm{\wti{z_1^j}}
\geq 1/960.\]

It follows that we may find a $(M_1,f)$-form $z_1^*$ in ${\bf Q}$ such that
 $z_1^*(\wti{z_1})=10^{-3}$
up to $1/k$.
 Furthermore there is an element $y_1^*$ in $B(\Yo)$, such that
\[\norm{y_1^*-z_1^*} \leq 1/f(M_1)\sum_{j=1}^{M_1}\Delta^{j}
\leq \Delta.10^{-3}/2 \leq \Delta \min\{\norm{y_1^*},\norm{z_1^*}\}.\]

Let $M_2=\sigma(z_1^*)$. As ${\bf Q}$ is dense, we had infinitely
many choices for $z_1^*$, so we may assume we chose $z_1^*$
such that 
 $1/2 f((M_2)^{1/40})^{1/2} \geq |ran(\wti{z_1})|$.
As in above, we may find a R.I.S. vector $z_2$ in $Z'$,
whose lifting $\wti{z_2}$ is of norm $1$, and a
$(M_2,f)$-form $z_2^*$ in ${\bf Q}$
such that $|z_2^*(z_2)-10^{-3}| \leq 1/k$, and such that
the couple $(\wti{z_2},z_2^*)$ is successive to $(\wti{z_1},z_1^*)$.
 Going on
in the same way, we build  sequences $(z_i)$,
$(\wti{z_i})$, $(z_i^*)$ and $(y_i^*)$ such that
for all $i$, $z_i$ is in $Z$ if $i$ is odd, in $Z'$ if $i$ is even,
$z_i^*$ is a $(M_i,f)$-form,
$\norm{y_i^*-z_i^*} \leq \Delta^{i}\min\{\norm{y_i^*},\norm{z_i^*}\}$
 $1/2 f((M_i)^{1/40})^{1/2} \geq |ran(\wti{z_{i-1}})|$, and
$(\wti{z_i},z_i^*)$ is successive to $(\wti{z_{i-1}},z_{i-1}^*)$.

 So by construction, $z_1^*,\ldots,z_k^*$ is a special sequence,
$(z_i^*)$ and $(y_i^*)$ are $\Delta$-equivalent (so that the sequence
$(z_i^*)$ is almost in $B(\Yo)$), and
$z_1,\ldots,z_k$ is a R.I.S. associated to 
$z_1^*,\ldots,z_k^*$.

It follows from Lemma 1 
 that $\norm{\sum_{i=1}^k y_i^*} \leq (1+2\Delta)\sqrt{f(k)}$
so

\[ \norm{\sum_{i=1}^k z_i} \geq
 \left((1+2\Delta)\sqrt{f(k)}\right)^{-1}
 \left(\sum_{i=1}^k y_i^*\right)\left(\sum_{i=1}^k z_i \right),\]

\[ \norm{\sum_{i=1}^k z_i} \geq
 \left((1+2\Delta)\sqrt{f(k)}\right)^{-1}
\left(\sum_{i=1}^k z_i^*(\wti{z_i})
-2\Delta \norm{\sum_{i=1}^k \wti{z_i}}\right),\]

\[ \norm{\sum_{i=1}^k z_i} \geq
 \left((1+2\Delta)\sqrt{f(k)}\right)^{-1} (10^{-3}k-1-2\Delta k) \geq 4.
 10^{-4}k/\sqrt{f(k)}.\]

On the other hand,
for all interval
 $E$, $|(\sum_{i=1}^k z_i^*)(\sum_{i=1}^k (-1)^{i}\wti{z_i})| \leq 2$,
so by Lemma 14, 
\[ \norm{\sum_{i=1}^k (-1)^{i} \wti{z_i}} \leq 7k/f(k).\]

It follows that 

\[ \norm{\sum_{i=1}^k (-1)^{i} z_i}
 \leq \norm{\sum_{i=1}^k (-1)^{i} \wti{z_i}}
+\sum_{i=1}^k \norm{z_i-\widehat{\wti{z_i}}}
 \leq 7k/f(k)+2\Delta \leq 8k/f(k).\]

If $z$ denotes the sum of the odd vectors, $z'$ the sum of
the even vectors, we have that $z \in Z$,
$z' \in Z'$, and
 \[\norm{z-z'} \leq 2.10^4
 f(k)^{-1/2} \norm{z+z'} \leq \delta \norm{z+z'}.\]
As 
$\delta$ is arbitrary, it follows that the sum of $Z$ and $Z'$
is not direct, and finally, that $X/Y$ is H.I..

\paragraph{Remark}
As $X$ is reflexive, it follows from Proposition 2 that
$X^*$ is H.I..

\section{There exists a H.I. space $\cal X$
such that ${\cal X}^*$ is not H.I.}
In this section, we build a H.I., not Q.H.I. space $\cal X$ as a certain sum
of two GM-type H.I. spaces $X_1$ and $X_2$.
 The space ${\cal X}^*$ will contain a direct sum of two subspaces.
By a simple generalization explained in the Appendix, it is even
possible to build a H.I. space ${\cal X}$ such that
${\cal X}^*$ contains a direct sum of $n$ subspaces.

\subsection{Proposition 3}
{\em For $i=1,2$, let $X_i$ be a H.I. Banach space, let $Z_i$ be
a subspace of $X_i$. Assume that $Z_1$ and $Z_2$ are isometric, and
that $X_1/Z_1$ and $X_2/Z_2$ are infinite dimensional and
totally incomparable. By abuse
of notation, we identify both $Z_1$ and $Z_2$ with a same space $Z$.
Let $\cal X$ be the quotient space $(X_1 \times X_2)/
\{(z,-z),z \in Z \}$. Then ${\cal X}$ is H.I. and
${\cal X}^*$ is not H.I..}
 
\underline{Proof} For $x_i$ in $X_i$, $i=1,2$, we denote by
$\hat{x_i}$ the class of $x_i$ in $X_i/Z_i$,
by  
$\widehat{(x_1,x_2)}$ the class of $(x_1,x_2)$ in $\cal X$. By definition,
\[ \norm{\widehat{(x_1,x_2)}}=\inf_{z \in Z}(\norm{x_1+z}+\norm{x_2-z}). \]
It follows that  the space ${\cal X}_1=\{\widehat{(x_1,0)},\ x_1 \in X_1\}$ is
isometric to $X_1$, the space
 ${\cal X}_2=\{\widehat{(0,x_2)},\ x_2 \in X_2\}$ is
isometric to $X_2$, and
the space ${\cal Z}=\{\widehat{(z,0)},\ z \in Z\}
=\{\widehat{(0,z)},\ z \in Z\}$ is
isometric to $Z$.
 As an easy consequence, we have the relation
\[ {\cal X}/{\cal Z} = {\cal X}_1/{\cal Z} \oplus {\cal X}_2/{\cal Z}
\simeq X_1/Z_1 \oplus X_2/Z_2,\]
so \[{\cal Z}^{\perp} \simeq (X_1/Z_1)^* \oplus (X_2/Z_2)^*,\]
and this proves that ${\cal X}^*$ is not H.I.. It remains
 to show
that $\cal X$ is H.I..
\
For $i=1,2$, we define a linear operator
 $\phi_i:{\cal X} \rightarrow X_i/Z_i$ by
$\phi_i(\widehat{(x_1,x_2)})=\widehat{x_i}$. It is easy to check that
$\phi_i$ is well defined. 
Now let $\cal W$ be a subspace of $\cal X$. There exists an $i$
such that $\phi_{i/{\cal W}}$ is infinitely singular: indeed, if
$\phi_{1/{\cal W}}$ and $\phi_{2/{\cal W}}$ are
both finitely singular, then there exists a subspace $\cal V$
of $\cal W$ on which $\phi_1$ and $\phi_2$ are isomorphisms into,
so that $X_1/Z_1$ and $X_2/Z_2$  have isomorphic subspaces,
a contradiction. 
\

Now assume for example that ${\phi_1}_{/{\cal W}}$ is infinitely
singular. As a consequence, there exists a norm $1$ sequence
$(w_n)_{n \in \N}$ such that
 $\phi_1(w_n) \stackrel{+\infty}{\rightarrow} 0$.
By definition of $\phi_1$, this means 
that $d(w_n,{\cal X}_2) \stackrel{+\infty}{\rightarrow} 0$.
 It follows easily
 that $\cal W$ and ${\cal X}_2$ have $Id+S$-isomorphic
subspaces. As ${\cal X}_2$ is isometric to $X_2$, it is H.I.;
it follows that $\cal Z$ is quasi-maximal in ${\cal X}_2$, so
$\cal W$ and $\cal Z$ have also
$Id+S$-isomorphic subspaces.
\

We have now proved that for every subspace $\cal W$ of $\cal X$,
$\cal W$ and $\cal Z$ have $Id+S$-isomorphic subspaces. This means 
that $\cal Z$ is quasi-maximal in $\cal X$. As $\cal Z$ is H.I.,
this implies that $\cal X$ is H.I..

\subsection{Proposition 4}
{\em For $i=1,2$, there exist $X_i$  reflexive Q.H.I. Banach space, 
$Z_i$ subspace of $X_i$, such that $Z_1$ and $Z_2$ are isometric, and
such that $X_1/Z_1$ and $X_2/Z_2$ are totally incomparable.}

\subsection{Definition of $X_1$ and $X_2$}
We shall define two spaces $X_i$ of the form $(c_{00},\norm{.}_i)$ for
$i=1,2$,
following a Gowers-Maurey's method, in which we force
$X_1$ and $X_2$ to have isometric subspaces $Z_1$ and $Z_2$.
The Banach spaces involved in  Proposition~4 are the completion
of the spaces $X_1$ and $X_2$.

Let $(e_n)_{n \in \N}$ be the canonical unit basis of $c_{00}$.

Let $\bf Q$ be the set of sequences with finite range,
rational coordinates  and maximum at most one in modulus.
We recall that $J$ is a set of integers $\{j_n,n \in \N\}$, such that
$f(j_1)>256$ and for all $n$, $\log \log \log j_{n+1} \geq 4 j_n^2$, that
$K=\{j_1,j_3,j_5,\ldots\}$, and
$L=\{j_2,j_4,j_6,\ldots\}$. Furthermore, we let
$L_1=\{j_2,j_6,j_{10},\ldots\}$,
$L_2=\{j_4,j_8,\ldots\}$.
For $i=1,2$, let $\sigma_i$ be an injection 
from the collection of finite sequences of successive
elements of ${\bf Q}$ to $L_i$.
For $i=1,2$, we choose $X_i$ to be of the form $(c_{00},\norm{.}_i)$,
and $Z_i=span\{e_{2n+1},n \in \N \}$. We may identify $Z_1$ and $Z_2$
with the same algebraic space $Z$, and $Z_1^{\perp}$ and $Z_2^{\perp}$ with
the same algebraic space 
$\Zo=span\{e_{2n}^*,n \in \N \}$. We now need some definitions.

\paragraph{Definitions}

A {\em dual couple} is a couple $(G,H)$ of
balanced bounded convex subsets of $c_{00}$.

Let $(G,H)$ be a dual couple.

A {\em $N$-Schlumprecht sum in $G$} is a
sum of the form $1/f(N) \sum_{i=1}^N x_i^*$, where the $x_i^*$'s
are successive in $G$. A {Schlumprecht sum in $G$} is a
 $N$-Schlumprecht sum in $G$ for some $N$.
 The set of Schlumprecht sums in $G$ is denoted
by $\Sigma(G)$.
In the same way, we define {\em Schlumprecht sums in $H$}.

A {\em special sequence in $G$}
is a sequence of successive vectors $x_1^*<\ldots<x_k^*$, with $k \in K$,
such that for $i=1,\ldots,k$, $x_{i}^*$ is an $M_i$-Schlumprecht sum in $G$
 with 
$M_i \geq j_{2k}$, and $M_i=\sigma_1(x_1^*,\ldots,x_{i-1}^*)$ for
$i=2,\ldots,k$.

A {\em special sum in $G$} is a sum of the form
$1/\sqrt{f(k)} \sum_{i=1}^k x_i^*$, where $x_1^* < \ldots < x_k^*$
is a special sequence in $G$. The set of special sums in $G$
is denoted by $S(G)$.

We similarly define {\em special sequences and special sums in $H$} 
replacing $\sigma_1$
by $\sigma_2$ in the above definition.

\

An {\em associated dual couple} is a dual couple $(G,H)$ such that
there exist two multivalued functions $a:G \rightarrow H$ and $b:H \rightarrow G$ satisfying
the four following properties.

\

(a) for all $x^* \in G$, all $y^* \in a(x^*)$, $y^*-x^*$ is in $\Zo$;
\

(b) for all $x^* \in G$, all $y^* \in a(x^*)$, $ran(y^*) \subset ran(x^*)$;
\

(c) for all $x^* \in G \cap \Zo$, $a(x^*)=\{0\}$;
\

(d) for all $N$-Schlumprecht sum $x^*$ in $G$, $a(x^*)$ contains an 
$N$-Schlumprecht sum in $H$,

\

and the four similar properties for $b$.

\paragraph{Definitions}

Let $(G,H)$ be an associated dual couple.

A {\em shadow sequence in $G$}
is a sequence of successive vectors $x_1^*<\ldots<x_k^*$ such that
there exists a special sequence $y_1^*<\ldots<y_k^*$
 in $H$
such that for all $i$, $x_i^*$ is an $M_i$-Schlumprecht sum in $G$ belonging
to $b(y_i^*)$, where
$M_i$ is the integer associated to $y_i^*$ in the definition of
the special sequence. A shadow sequence in $G$ can
be thought of as the "shadow" of a special sequence in $H$. 

A {\em shadow sum in $G$} is a sum of the form
$1/\sqrt{f(k)} \sum_{i=1}^k x_i^*$, where $x_1^* < \ldots < x_k^*$
is a shadow sequence in $G$. The set of shadow sums in $G$
is denoted by $s(G,H)$.

We similarly define {\em shadow sequences and shadow sums in $H$}, and
denote the set of shadow sums by $s(H,G)$.

\

 To define the norms,
we shall now build by induction an associated dual couple $(C,D)$ 
where $C$ (resp. $D$) is meant
 to be almost the dual unit
ball of $X_1$ (resp. $X_2$).
 We shall build $C$ as $\cup_{n \in \N}C_n$, building the
increasing sequence $C_n$ by induction. We shall also build
$a$ by induction, defining a function $a_n$ from $C_n$ to $D_n$
at each step $n$;
but to simplify the notation, we shall denote all the terms
of the sequence by $a$ (and we shall do symmetrically the same 
for $D$ and $b$).

In this situation, Property (a) ensures
 that the subspaces $Z_1$ and $Z_2$ are isometric.
Properties (b) and (d) allow us to give convenient
properties to the images
 by $a$ of the special sequences,
that is the shadow sequences. Property (c) introduces an asymmetry,
and as a consequence, the quotient spaces $X_1/Z_1$ and
$X_2/Z_2$ will be totally incomparable.

\paragraph{Construction}

At the first step, we define $(C_0,D_0)$ 
to be $(B(l_1),B(l_1))$,
 $a$ and $b$ by 
$a(\sum_{i \in \N} \la_i e_i^*)=b(\sum_{i \in \N} \la_i e_i^*)=
\sum_{i\ odd} \la_i e_i^*$.
It is easy to check that $(C_0,D_0)$ is an associated dual couple.

\

Now assume we are given an associated dual couple $(C_{n-1},D_{n-1})$, with
 functions
$a:C_{n-1} \rightarrow D_{n-1}$ and $b:D_{n-1} \rightarrow C_{n-1}$.

We define $C^{\prime}_{n-1}$ to be
$\Sigma(C_{n-1}) \cup S(C_{n-1}) \cup s(C_{n-1},D_{n-1})$, and
$C_n$ to be the set of elements of the 
form $E(\sum_{i=1}^M \la_i x^*_i)$, where $E$ is an interval projection,
$\sum_{i=1}^M |\la_i|=1$, and
 for all $i$, $x^*_i$ is in $C^{\prime}_{n-1}$.
We define $D_n$ in a similar way.

We now extend $a$ to $C_n$. Let $x^* \in C_n$.
If $x^* \in \Zo$, then we let $a(x^*)=\{0\}$. We now define a construction
if $x^*$ is not in $\Zo$.

The set $a(x^*)$ may be already defined or not (it is when
$x^*$ is in $C_{n-1}$); if not we may assume $a(x^*)=\emptyset$. Then
we add new values to the set $a(x^*)$ in each of the
following cases (notice that at least one of the possibilities happens, so
that $a$ is well defined on the whole of $C_n$, but that the possibilities
are not exclusive).

- If $x^*$ is a Schlumprecht sum of
 the form $f(N)^{-1} \sum_{i=1}^N x_i^*$ with $x_i^* \in C_{n-1}$
then we add to $a(x^*)$  
the set $f(N)^{-1} \sum_{i=1}^N a(x_i^*)$.

- If $x^*$ is a special sum of 
the form $f(k)^{-1/2} \sum_{i=1}^k x_i^*$
 where $x_i^*$ is an $(M_i,f)$-form in $C_{n-1}$ then we add to the
set $a(x^*)$
the set of all sums of the form $f(k)^{-1/2} \sum_{i=1}^k y_i^*$,
whe\-re $y_i^*$ is an $(M_i,f)$-form in $a(x_i^*)$.

- If $x^*$ is a shadow sum of 
the form $f(k)^{-1/2} \sum_{i=1}^k x_i^*$ with
$x_i^* \in b(y_i^*)$ and $y_1^*,\ldots,y_k^*$ is a special
sum in $D_{n-1}$,
 then we add to the set $a(x^*)$
the singleton $\{ f(k)^{-1/2} \sum_{i=1}^k y_i^*\}$.      

- If $x^*$ is the projection of a convex combination
of elements of the three previous forms, that is,
$x^*=ran(x^*)(\sum_i \la_i x_i^*)$,
then we add to the set $a(x^*)$ the set $ran(x^*)(\sum_i \la_i a(x_i^*))$,
 $a(x_i^*)$ being
defined as above whether $x_i^*$ is a Schlumprecht sum,
a special sum, or a shadow sum in $C_{n-1}$. Notice that 
we only use this construction when $x^*$ is not in $Z^{\perp}$.

It is then easy to check that $a$ (resp. $b$) takes its values in 
$D_n$ (resp. in $C_n$) and
that it still satisfies the four properties (a)-(d), so $(C_n,D_n)$ is a
dual couple. 

\

We finally define $C$ as $\cup_{n \in \N}C_n$
and $D$ as $\cup_{n \in \N}D_n$; the multifunction
$a$ (resp. $b$) is defined on $C$ (resp. $D$), so
$(C,D)$ is a dual couple.

\

We define $X_1$ by its norm $\norm{.}_1=\sup_{x^* \in C}<x^*,.>$,
$X_2$ by its norm $\norm{.}_2=\sup_{y^* \in D}<y^*,.>$.

\subsection{$Z_1$ and $Z_2$ are isometric subspaces}

Let $z$ be an element of $Z$. Then
\[ \norm{z}_1=\sup_{x^* \in C}(<x^*,z>) \leq
\sup_{x^* \in C,y^* \in a(x^*)}(<x^*-y^*,z>+<y^*,z>).\]
Now by definition of $a$, $x^*-y^* \in \Zo$, so
$<x^*-y^*,z>=0$; and as $y^*$ is in $D$, $(<y^*,z>) \leq \norm{z}_2$.
It follows that $\norm{z}_1 \leq \norm{z}_2$, and by symmetry,
$\norm{z}_1=\norm{z}_2$.

\subsection{$X_1$ and $X_2$ are reflexive and Q.H.I.}

By symmetry, it is enough to show that for example $X_1$ is reflexive 
and Q.H.I..

\paragraph{Remarks}

With the definition following Lemma 7, a special sequence 
in $X_1$ (resp. $X_2$) is an $L_1$ (resp. $L_2$)-sequence,
a shadow sequence in $X_1$ (resp. $X_2$) is an
 $L_2$ (resp. $L_1$)-sequence.
 Notice also that 
the space $X_1$ (resp. $X_2$) has GM-type,
the set ${\cal S}_1$ (resp. ${\cal S}_2$) 
being the set of special and shadow sums, and so it is reflexive.

\paragraph{Lemma 15}
{\em Let $X$ have GM-type. Let $x_1^*,\ldots,x_k^*$
be a $L_1$-sequence in $X^*$.
Let $x_1<\ldots<x_k$ in $X$ be a R.I.S. associated to $x_1^*,\ldots,x_k^*$.
Let $x=x_1+\ldots+x_k$.
Then for every $L_2$-sum $z^*$ of length $k$ in $X^*$,
 every interval $E$,
$|z^*(Ex)| \leq 1/4$.}

\underline{Proof}

Let $z^*$ be an $L_2$-sum, $E$ be an in\-ter\-val.
Then $z^*=1/\sqrt{f(k)} \sum_{i=1}^k z_i^*$, whe\-re $z_i^*$
is a $(l_i,f)$-form, and $l_i$ is in $L_2$.
 For eve\-ry $j$, $x_j$ is a norm $1$ R.I.S.
of leng\-th in $L_1$. As $L_1$ and $L_2$ are dis\-joint,
and the lengths are greater than $j_{2k}$ by definition of an
$L$-sequence and of an
associated R.I.S.,  it fol\-lows
from Lem\-ma 11 that $|z_i^*(Ex_j)| \leq 1/k^2$.
Fi\-nal\-ly,
\[ |z^*(Ex)| \leq 1/\sqrt{f(k)} \leq 1/4.\]

\paragraph{Remark}
The result is also true for an $L_1$-sum and a R.I.S. associated to
a $L_2$-sequence.

\paragraph{Lemma 16}
{\em Let $x_1^*,\ldots,x_k^*$ be a special sequence in $X_1^*$.
Let $x_1<\ldots<x_k$ be associated to $x_1^*,\ldots,x_k^*$.
Let $x=\sum_{i=1}^k x_i$.

Assume that for
every interval $E$,
 $|(\sum_{i=1}^k x_i^*)(Ex)| \leq 2$, then
\[ \norm{x} \leq 7k/f(k). \]}

\underline{Proof} 
We already know that $X_1$ has GM-type, so from Lemma 15, 
for every shadow sum $z^*$,
 every interval 
$E$, $|z^*(Ex)| \leq 1/4$.
We may now apply the same proof as in Lemma 14, replacing
$\sigma$ by $\sigma_1$.
 Lemma 12 then gives the conclusion.

\paragraph{Proposition}
{\em $X_1$ is Q.H.I..}

\underline{Proof} As we said above, $X_1$ is of GM-type, like
Gowers-Maurey's space. Furthermore, we have proved Lemma 16,
which is similar to Lemma 14 in the case of Gowers-Maurey's space.
It follows that our final proof that Gowers-Maurey's space
is Q.H.I. is still valid in the case of $X_1$.

\subsection{$X_1/Z_1$ and $X_2/Z_2$ are totally incomparable.}

\paragraph{Lemma 17}
{\em Let $y_1^*,\ldots,y_k^*$ be a special sequence in 
$Z_2^{\perp}$.
Let $x_1<\ldots<x_k$ in $X_1$ be associated to
 $y_1^*,\ldots,y_k^*$.
Let $x=\sum_{i=1}^k x_i$. Then
\[ \norm{x} \leq 7 k/f(k). \]}

\underline{Proof} The space $X_1$ is of GM-type, so,
 because of Lemma 12, it
 is enough to show
that for every interval $E$, every special or shadow sum
$z^*$ in $X_1^*$, $|z^*(Ex)| \leq 1/4$. By Lemma 15, this is true
for every special sum. Now let $E$ be an interval and
let $z^*$ be a shadow sum in $X_1^*$.

 For every $i$, let $M_i$ be
such that $y_i^*$ is an $(M_i,f)$-form.
There exists a special sequence $t_1^*,\ldots,t_k^*$ in $X_2^*$
such that $z^*=1/\sqrt{f(k)} \sum_{i=1}^k z_i^*$ with
for every $i$, $z_i^* \in b(t_i^*)$;
let $N_i$ be
such that $t_i^*$ is an $(N_i,f)$-form; by construction
$z_i^*$ is also an $(N_i,f)$-form.
Let $I=sup\{i/ M_i=N_i\}$, or $0$ if no such $I$ exists.
For $i<I$, because $\sigma_2$ is an injection, we have
that $t_i^*=y_i^*$. It follows that $t_i^*$ is in $Z_2^{\perp}$, so
$b(t_i^*)=\{0\}$, and $z_i^*=0$.
For $i>I$, $z_i^*$ is an $(N_i,f)$-form, and $N_i$ is an element
of $L_2$ greater than $j_{2k}$ and different from all $M_j$, $j=1,\ldots,k$.
It follows by Lemma 11 that $|z_i^*(Ex_j)|\leq 1/k^2$.
Finally,
\[|z^*(Ex)| \leq 1/\sqrt{f(k)} (0+|z_I^*(x_I)|
+k^2.k^{-2}) \leq 2/\sqrt{f(k)} \leq 1/4.\]

\paragraph{Total incomparability}

We now assume that there exists an isomorphism $\al$ between
a subspace $W_1$ of $X_1/Z_1$ and a subspace $W_2$ 
of $X_2/Z_2$ and we intend
to find a contradiction.

First we notice that
$X_2/Z_2$ has a basis (namely the basis $(e'_{2n})_{n \in \N}$ 
dual to the basis
$(e_{2n}^*)_{n \in \N}$). Let $(w_n)$ be an average basic sequence in $W_1$.  Up to a perturbation on $\alpha$, we
may assume that the sequence $(\alpha(w_n))$ is a  sequence
of unit vectors, successive with respect to $(e'_{2n})$.

Now let $k \in K$. In $span\{(w_n)\}$, we
 may find a R.I.S vector $x_1=\sum_{i=1}^{M_1} x_1^{i}$ in $X_1$ 
of length $M_1=j_{2k}$, such that
$\wti{x_1}=\sum_{i=1}^{M_1} \wti{x_1^{i}}$ is of norm $1$.
For $i=1,\ldots,M_1$, let $y_1^{i*} \in Z_2^{\perp}$ be a
functional that norms $\alpha(x_1^{i})$ and such that $ran(y_1^{i*})
\subset ran(\alpha(x_1^{i}))$, and let
$y^{\prime *}_1$ be the $(M_1,f)$-form
$f(M_1)^{-1} \sum_{i=1}^{M_1} y_1^{i*}$.
We have that $y_1^{\prime*}$ is in $Z_2^{\perp}$ and
\[y_1^{\prime*}(\alpha(x_1))=f(M_1)^{-1}\sum_{i=1}^{M_1}\norm{\alpha(x_1^{i})},\]
so by Lemma 7 and Lemma 9,
\[y_1^{\prime*}(\alpha(x_1)) \geq
 (4\norm{\alpha^{-1}}f(M_1))^{-1}\sum_{i=1}^{M_1}\norm{\wti{x_1^{i}}}
\geq (20\norm{\alpha^{-1}})^{-1}.\]

 We
 may then choose a $(M_{1},f)$-form
$y^*_{1}$ in $B(Z_2^{\perp})$ such that
$y^*_{1}$ is in ${\bf Q}$, such that $ran(y^*_{1}) \subset
ran(\alpha(x_{1}))$, and such that
$y^*_{1}(\alpha(x_{1}))\geq (25\norm{\alpha^{-1}})^{-1}$.

 We then define $M_2=\sigma_2(y^*_1)$, and 
 we may assume we chose $y^*_{1}$ such that
$1/2 f((M_2)^{1/40})^{1/2} \geq |ran(\wti{x_{1}})|$.
We now let
$x_2$  be a R.I.S. of length $M_2$ in $span\{(w_n)\}$, such that
 $ran(\wti{x_2}) > ran(\wti{x_1})$.
Following in the same way, we obtain sequences
$x_i$, $y^*_i$ for $i=1,\ldots,k$, such that
the sequence $y^*_1,\ldots,y^*_k$ is
a special sequence in $Z_2^{\perp}$ and
$\wti{x_1},\ldots,\wti{x_k}$ is associated to $y^*_1,\ldots,y^*_k$.

It follows that
\[\norm{\alpha(\sum_{i=1}^k x_i)}
 \geq f(k)^{-1/2} \sum_{i=1}^k y_i^*(\alpha(x_i))
 \geq (25\norm{\alpha^{-1}})^{-1}kf(k)^{-1/2},\]
while
as $(\wti{x_i})$ is associated to $(y_i^*)$, by Lemma 17, 
\[ \norm{\sum_{i=1}^k \wti{x_i}} \leq 7 kf(k)^{-1},\]
and
\[ \norm{\sum_{i=1}^k x_i} \leq 7 kf(k)^{-1}+2\Delta \leq 8kf(k)^{-1},\]

It follows that 
$\norm{\alpha} \geq (200\norm{\alpha^{-1}})^{-1}\sqrt{f(k)}$, and 
this for any $k$,
contradicting the boundedness of the operator $\alpha$.

\section{Appendix}
We give a sketch of the proof of the existence of a H.I. space
$\cal X$ such that ${\cal X}^*$ contains a direct sum of $n$
subspaces.

\subsection{Proposition}
{\em For $i=1,\ldots,n$, let $X_i$ be a H.I. Banach space, let $Z_i$ be
a subspace of $X_i$. Assume that the spaces $Z_i$ are all isometric
to a same space $Z$, and
that for any $i \neq j$, $X_i/Z_i$ and $X_j/Z_j$ are infinite dimensional and
totally incomparable.
Let $Z_{[1,n]}=\{ (z_1,\ldots,z_n) \in Z_1 \times \cdots \times Z_n
\ /\ \sum_{i=1}^n z_i=0\}$.
Let $\cal X$ be the quotient space $(X_1 \times \cdots \times X_n)/
Z_{[1,n]}$. Then ${\cal X}$ is H.I. and
${\cal X}^*$ contains a direct sum of $n$ subspaces.}
 
\underline{Proof} 

We use the same notation as in the case $n=2$, in particular
we let ${\cal Z}=\{\widehat{(z,\ldots,0)},\ z \in Z\}$, and we show 
that 
 \[{\cal Z}^{\perp} \simeq \bigoplus_{i=1}^n (X_i/Z_i)^*.\]

Now we 
consider $\cal W$ a subspace of $\cal X$. There exists at most
one value $i_W$ of $i$ such that $\phi_{i/{\cal W}}$ is finitely singular,
otherwise two quotient spaces $X_i/Z_i$ and $X_j/Z_j$ would have 
isomorphic subspaces.

 It follows easily
that $\cal W$ and ${\cal X}_{i_W}$ have $Id+S$-isomorphic
subspaces, and finally that ${\cal X}$ is H.I..

\subsection{Proposition}
{\em For $i=1,\ldots,n$, there exist $X_i$  Q.H.I. reflexive Banach space, 
$Z_i$ subspace of $X_i$, such that all $Z_i$ are isometric, and
such that for any $i \neq j$,
 $X_i/Z_i$ and $X_j/Z_j$ are totally incomparable.}

We make a construction similar to the case $n=2$, using a partition
of $L$ in $n$ subsets $L_1,\ldots,L_n$.
We build $n$ balanced bounded convex subsets $C_1,\ldots,C_n$ of $c_{00}$,
 and multifunctions $a_{ij}:C_i \rightarrow C_j$ for $i \neq j$,
such that for all $i \neq j$, $(C_i,C_j)$ is an associated dual couple.
The difference is that we have $n-1$ kinds of shadow sequences
in each $C_i$ (those coming  special sequences in $C_j$ for all
$j \neq i$). Defining $a_{ij}$ for the new kinds of shadow sums
is not difficult. The proof of the proposition then follows exactly
the case $n=2$.

\subsection{Remark}
We recall a definition from \cite{F}: a Banach space
is said to be $HD_n$ if the maximum number of subspaces
in a direct sum is finite and equal to $n$.
Here, for $i=1,\ldots,n$, $X_i$ is reflexive Q.H.I., so
by Proposition 2, $X_i^*$ is H.I., so
 it follows from Corollary 1 of \cite{F} that
$X_1^* \oplus \ldots \oplus X_n^*$ is $HD_n$.
Now ${\cal X}^*$ is a subspace of $X_1^* \oplus \ldots \oplus X_n^*$,
so it is $HD_m$ for some $m \leq n$,
and it contains the direct sum of $n$ subspaces ${\cal Z}^{\perp}$, so
$m=n$; but ${\cal X}^*$ is not decomposable, otherwise ${\cal X}$
would be decomposable. So ${\cal X}^*$ is an example
of a non decomposable $HD_n$ space.

\end{document}